\documentclass[final,12pt]{elsarticle}

\usepackage{amssymb,amsmath,algpseudocode}
\newtheorem{theorem}{Theorem}
\newtheorem{lemma}[theorem]{Lemma}
\newtheorem{corollary}[theorem]{Corollary}
\usepackage{epstopdf}

\journal{Linear Algebra And Its Applications}

\begin{document}

\begin{frontmatter}


%%s

\title{Approximation of the Scattering Amplitude using Nonsymmetric Saddle Point Matrices}


\author[addrusm]{Amber S. Robertson}
\author[addrusm]{James V. Lambers}

\address[addrusm]{Department of Mathematics, University of Southern Mississippi, 118 College Dr \#5045, Hattiesburg, MS 39406, USA}

\begin{abstract}
In this paper we examine iterative methods for solving the forward
($A{\bf x}={\bf b}$)
and adjoint ($A^{T}{\bf y}={\bf g}$)
systems of linear equations used to approximate the scattering amplitude, defined by ${\bf g}^{T}{\bf x}={\bf y}^{T}{\bf b}$.  Based on an idea first proposed by Gene Golub, we use a conjugate gradient-like iteration for a nonsymmetric saddle point matrix that is constructed so as to have a real positive spectrum. Numerical experiments show that this method is more consistent than known methods for computing the scattering amplitude such as GLSQR or QMR. 
%Then, we use techniques from "matrices, moments, and quadrature" to compute the scattering amplitude without solving the system directly.
We then demonstrate that when combined with known preconditioning techniques, the proposed method exhibits more rapid convergence than state-of-the-art iterative methods for nonsymmetric systems.
\end{abstract}


\end{frontmatter}

% \linenumbers

\section{INTRODUCTION}
\label{sec1}

\subsection{The Scattering Amplitude Problem}
The core objective of this paper is to design and implement an iterative method for the solution of a system where the coefficient matrix 
%is a nonsymmetric saddle point matrix, meaning it 
is large, sparse, and nonsymmetric. The proposed method should be more efficient and robust than existing methods for solving such systems.
One application in which such a system arises is in the computation of the \emph{scattering amplitude}.
The scattering amplitude, in quantum physics, is the amplitude of the outgoing spherical wave relative to that of the incoming plane wave \cite{GSW}. It is useful when it is of interest to know what is reflected when a radar wave is impinging on a certain object. The scattering amplitude can be computed by taking the inner product of the right hand side vector ${\bf g}$ of the $\emph{adjoint system}$
\begin{equation}
A^{T}{\bf y}={\bf g}
\end{equation}
and the solution {\bf x} of the $\emph{forward system}$
\begin{equation}
\label{ax=b}
A{\bf x}={\bf b}.  
\end{equation}
Applications of the scattering amplitude come up in nuclear physics \cite{ARN}, quantum mechanics \cite{LALI}, and computational fluid dynamics (CFD) \cite{GIPI}. One particular application is in the design of stealth planes \cite{ARN}.

The scattering amplitude ${\bf g}^{T}{\bf x}={\bf y}^{T}{\bf b}$ creates a relationship between the right hand side of the adjoint system and the solution to the forward system in signal processing. The field ${\bf x}$ is determined from the signal ${\bf b}$ in the system $A{\bf x}={\bf b}$. Then the signal is received on an antenna characterized by the vector ${\bf g}$ which is the right hand side of the adjoint system $A^{T}{\bf y}={\bf g}$, and it is expressed as ${\bf g}^{T}{\bf x}$ \cite{GSW}. We are interested in efficiently approximating the scattering amplitude. It is informative to look at methods that other researchers have used to solve this problem, which will be discussed below. 

The solution of the linear system (\ref{ax=b}) is important for many applications beyond the scattering amplitude, such as in the numerical solution of PDE with non-self-adjoint spatial differential operators. This solution can be obtained in many different ways, depending on the properties of the matrix $A$. The $LDL^{T}$ factorization can be used to solve some problems with a symmetric matrix or a \emph{Cholesky factorization} can be used if the matrix is also known to be positive definite \cite{GVL}. However, for large, sparse systems, an iterative method is preferred. The \emph{conjugate gradient} method is the preferred iterative method for a symmetric positive definite matrix $A$ \cite{GVL}. However it is much more difficult to find this solution for a matrix that is not symmetric positive definite. In the case that we have a matrix that is not symmetric, we can use methods like the \emph{biconjugate gradient (BiCG)} \cite{BRRZ} and \emph{generalized minimal residual (GMRES)} methods \cite{GMRES}. If we have a matrix that is symmetric but indefinite, \emph{SymmLQ} \cite{SYMLQ,SYMLQ2} is the iterative method of choice. Since the scattering amplitude depends on both the forward and adjoint problem, it makse sense to use methods that take both the forward and adjoint problems into account, like the \emph{quasi-minimal residual (QMR)} \cite{QMR} and \emph{generalized least squares residual (GLSQR)} methods\cite{GLSQR}.\\

\subsection{Approximation of the Scattering Amplitude}
%\indent Another important objective of this paper is to develop methods that approximate expressions of bilinear form,
%$${\bf u}f(A){\bf v}$$ where {\bf u} and {\bf v} are $n$ vectors and $A$ is an $n\times n$ nonsymmetric matrix. The scattering amplitude can be expressed as a bilinear form where $f(A)$ represents the inverse operator. 
\indent The method of this paper employs a conjugate gradient-like approach since, for large, sparse matrices, it is best to use an iterative approach, such as the conjugate gradient method
\cite{CGHS} which is particularly effective for symmetric positive definite matrices. In particular, conjugate gradient has a very rapid convergence if $A$ is near the identity either in the sense of a low rank perturbation or in the sense of the norm. In \cite{GVL} it is stated that
\begin{theorem}
 If $A=I + B$ is an $n \times n$ symmetric positive definite matrix and rank($B$)=$ r$ then the Hestenes-Stiefel conjugate gradient algorithm converges in at most $r+1$ steps.
\end{theorem}
\begin{theorem}
Suppose $ A \in \mathbb{R}^{n \times n}$ is symmetric positive definite and $b \in \mathbb{R}$. If the Hestenes-Stiefel algorithm produces iterates ${\bf x}_{k}$ and $\kappa =\kappa_{2}(A)$ then
$$
\|{\bf x}-{\bf x}_{k}\|_{A} \leq 2 \|{\bf x}-{\bf x}_0 \|_{A} \left( \frac{ \sqrt{ \kappa}-1}{\sqrt {\kappa}+1} \right) ^{k},
$$
where $\|{\bf w}\|_{A}=\sqrt{{\bf w}^{T}A{\bf w}}$.
\end{theorem}
It is also stated in \cite{GVL} that the accuracy of ${\bf x}_{k}$ is often better than this theorem predicts and that the conjugate gradient method converges very rapidly in the $A$-norm if  $\kappa_{2}(A) \approx 1$, where $\kappa_{2}(A)$ is the $\emph{condition number}$ of $A$, defined by
$$\kappa_{2}(A)=\|A\|_{2}\|A^{-1}\|_{2}=\frac{\sigma_{\max}(A)}{\sigma_{\min}(A)}$$ 
with $\sigma_{\max}$ and $\sigma_{\min}$ referring to the largest and smallest singular values, respectively.

Multiplying both sides of $A{\bf x} = {\bf b}$ by $A^T$ yields the normal equations with a symmetric matrix $A^T A$ that is also positive definite when $A$ is invertible. 
However, this approach is not conducive to solving the forward and adjoint problems simultaneously.
Furthermore, a significant problem with using $A^{T}A$ is that now the condition number in the two-norm 
%that is also equal to the largest singular value over the smallest, 
is squared for $A^{T}A$. Since this increases the sensitivity of the matrix, possibly making it ill-conditioned, this paper explores an alternative approach. The idea is to transform the problems $A{\bf x} = {\bf b}$ and $A^T {\bf y} = {\bf g}$ into an equivalent system in which the matrix can be guaranteed to have real, positive eigenvalues, as well as eigenvectors that are in some sense orthogonal, which is then conducive to solution using a conjugate gradient-like iteration. It is not necessarily symmetry that we seek, but we will have symmetry with respect to some inner product. 
To this end, we use an idea first proposed by Gene Golub in \cite{JLGG}, and consider a nonsymmetric saddle point matrix that has the form
$$
M=\left[ \begin{array}{cc}
A^{T}WA & A^{T}  \\
-A & 0 \\
\end{array} \right].
$$

As required by the definition of a nonsymmetric saddle point matrix, 
we assume that the matrix $W$ is symmetric positive definite.  The goal is to choose $W$ so that we can guarantee $M$ has real, positive eigenvalues. 
In this paper we will introduce the \emph{nonsymmetric saddle point conjugate gradient} (NspCG) method to solve a nonsymmetric, large, sparse linear system, which will then allow us to compute the scattering amplitude. We will also use ILU preconditioning with NspCG, which gives rapid convergence compared to existing methods for solving such systems.

%In this paper we use several different approaches to compute the scattering amplitude: nonsymmetric saddle point conjugate gradient, unsymmetric Lanczos with perturbation of the initial vectors, and unsymmetric block Lanczos. These iterative approaches to solving the linear system and computing the scattering amplitude are used with the matrix $M$. We see some improvement in the rate of convergence for some of these methods. 
%We will also try some of these approaches, (unsymmetric Lanczos, block GLSQR, and symmetrizing of the initial vectors) with a symmetric matrix $C$, where
%$$
%C=
%\left[ \begin{array}{cc}
%0 & A \\
%A^{T} & 0\\
%\end{array} \right].
%$$
%The scattering amplitude is ${\bf u }^{T}C^{-1}{\bf v}$ where, 
%$$
%{\bf u}=
%\left[ \begin{array}{c}
%{\bf b}\\
%{\bf 0}\\
%\end{array} \right]
%\quad \quad 
%{\bf v}=
%\left[ \begin{array}{c}
%{\bf 0}\\
%{\bf g}\\
%\end{array} \right].
%$$
\indent This paper is organized as follows. In Section \ref{sec2} we discuss the known methods for solving a large linear system with iterative approaches to compute the scattering amplitude such as Bidiagonalization or least squares QR (LSQR), quasi minimum residual (QMR), and block generalized LSQR (GLSQR). 
%Section \ref{sec3} will give some necessary background on matrices, moments, and quadrature (MMQ). The MMQ background is necessary to understand the methods of this paper because we want to approximate an expression of the form
%$$
%{\bf u}^{T}f(A){\bf v},
%$$
%where ${\bf u}$ and ${\bf v}$ are $N$-vectors, $f(\lambda)=\lambda^{-1}$ , and $A$ is an $N \times N$ matrix.
 In Section \ref{sec4} we will introduce the method of this paper, NspCG.
Section 4 will include an analysis of the numerical results. The preconditioning techniques and results can be found in Section \ref{sec5b}. The conclusions and discussion of possible future work will be given in Section \ref{sec6}.

%THIS IS FILE diss_chap2.tex
%
%THIS IS THE NEXT CHAPTER. 
%IT ILLUSTRATES THE USE OF SECTION HEADERS AND EQUATIONS.

\section{Methods for Solving the Linear Systems of the Forward and Adjoint Problems} 
\label{sec2}

\subsection{QMR approach}
The QMR approach \cite{QMR,GSW} is based on the spectral decomposition $ A=XDX^{-1}$; also
the basis of the QMR approach is the unsymmetric Lanczos \cite{GVL,UnsymmL} process which generates two sequences
$$ 
V_{k}=
\left[ \begin{array}{cccc}
{\bf v}_{1} & {\bf v}_{2} & \ldots & {\bf v}_{k}
\end{array} \right]
$$ 
$$ 
W_{k}=
\left[ \begin{array}{cccc}
{\bf w}_{1} & {\bf w}_{2} & \ldots & {\bf w}_{k}
\end{array} \right]
$$ 
that are biorthogonal, meaning $V_{k}^{T}W_{K}=I$. We have the following relations:
\begin{eqnarray}
 AV_{k}&=&V_{k+1}T_{k+1,k},\\
A^{T}W_{k}&=&W_{k+1}\hat{T}_{k+1,k}.
\end{eqnarray}
where the tridiagonal matrices 
$$T_{k+1,k}=
\left[ \begin{array}{ccccc}
\alpha_{1} & \gamma_{1} & \, & \, & \, \\
\beta_{1} & \alpha_{2} & \gamma_{2} &\, & \, \\
\, & \beta_{2} & \ddots & \ddots & \, \\
\, & \, & \ddots & \ddots & \gamma_{k-1} \\
\, & \, & \, & \beta_{k-1} & \alpha_{k} \\
\, & \, & \, & \, & \beta_{k} \\
\end{array} \right]=
\left[ \begin{array}{c}
T_{k,k}\\
\beta_{k}{\bf e}_{k}^{T}\\
\end{array} \right]
$$ 
and 
$$\hat{T}_{k+1,k}=
\left[ \begin{array}{ccccc}
\hat{\alpha}_{1} & \hat{\gamma}_{1} & \, & \, & \, \\
\hat{\beta}_{1} & \hat{\alpha}_{2} & \hat{\gamma}_{2} &\, & \, \\
\, & \hat{\beta}_{2} & \ddots & \ddots & \, \\
\, & \, & \ddots & \ddots & \hat{\gamma}_{k-1} \\
\, & \, & \, & \hat{\beta}_{k-1} & \hat{\alpha}_{k} \\
\, & \, & \, & \, & \hat{\beta}_{k} \\
\end{array} \right]=
\left[ \begin{array}{c}
\hat{T}_{k,k}\\
\hat{\beta}_{k}{\bf e}_{k}^{T}\\
\end{array} \right]
$$ have block structures in which $T_{k,k}$ and $\hat{T}_{k,k}$ are not necessarily symmetric.

The residual, ${\bf r}={\bf b}-A{\bf x}$, in each iteration can be expressed as 
\begin{eqnarray}
\nonumber
\|{\bf r}_{k}\|&=&\|{\bf b}-A{\bf x}_{k}\| \\
\nonumber
&=&\|{\bf b}-A{\bf x}_{0}-AV_{k}{\bf c}_{k}\| \\
\nonumber
&=&\| {\bf r}_{0}-V_{k+1}T_{k+1,k}{\bf c}_{k}\| \\
&=&\|V_{k+1}(\|{\bf r}_{0}\| {\bf e}_{1}-T_{k+1,k}{\bf c}_{k})\| 
\end{eqnarray} 
with a choice of ${\bf v}_{1}=\frac{{\bf r}_{0}}{\|{\bf r}_{0}\|}$ where ${\bf r}_{0}={\bf b}-A{\bf x}_{0}$ and ${\bf x}_{k}={\bf x}_{0}+V_{k}{\bf c}_{k}$. 
We now have the quasi-residual $\|{\bf r}_{k}^{Q}\|= \|\|{\bf r}_{0}\|{\bf e}_{1}-T_{k+1,k}{\bf c}_{k}\|.$ Then we choose ${\bf w}_{1}=\frac{\bf{s}_{0}}{\|{\bf s}_{0}\|}$, where ${\bf s}_{0}={\bf g}-A^{T}{\bf y}_{0}$ and ${\bf y}_{k}={\bf y}_{0}+{\bf w}_{k}{\bf d}_{k}$. Then the adjoint residual is 
$\|{\bf s}_{k}^{Q}\|=\| \| {\bf s}_{0}\|{\bf e}_{1}-\hat{T}_{k+1,k}{\bf d}_{k} \|$.
The vectors ${\bf c}_{k}$ and ${\bf d}_{k}$ are the solutions of the least squares problems for minimizing $\|{\bf r}_{k}^{Q}\|$ and $\|{\bf s}_{k}^{Q}\|$. 
So now the solutions can be defined as 
\begin{eqnarray}
{\bf x}_{k}&=&{\bf x}_{0}+V_{k}{\bf c}_{k}\\
{\bf y}_{k}&=&{\bf y}_{0}+U_{k}{\bf d}_{k}.
\end{eqnarray}

\subsection{LSQR approach}
In LSQR \cite{GSW,SYMLQ2}, a truncated bidiagonalization is used in order to solve the forward and adjoint problems approximately. The bidiagonal factorization of $A$ is given by $A=UBV^{T}$ where $U$ and $V$ are orthogonal and $B$ is bidiagonal. Thus the forward and adjoint systems can be written as 
\begin{equation}
\label{2.1}
UBV^{T}{\bf x}= {\bf b} 
\end{equation}
\begin{equation}
\label{2.2}
VB^{T}U^{T}{\bf y}={\bf g}.
\end{equation}
Now we can solve (\ref{2.1}) by solving the following two systems
\begin{eqnarray}
B{\bf z}&=&U^{T}{\bf b}\\
{\bf x}&=&V^{T}{\bf z},
\end{eqnarray}
and we can solve (\ref{2.2}) by solving
\begin{eqnarray}
B{^{T}\bf w}&=&V^{T}{\bf g}\\
{\bf y}&=&U^{T}{\bf w}.
\end{eqnarray}
%Solving $$A{\bf x}={\bf b} \quad\quad A^{T}{\bf y}={\bf g}$$ simultaneously can be done by %solving  
%$$ \left[ \begin{array}{cc}
%0 & A  \\
%A^{T} & 0 \\
%\end{array} \right]
 %\left[ \begin{array}{cc}
%{\bf y}  \\
%{\bf x} \\
%\end{array} \right]=
 %\left[ \begin{array}{cc}
%{\bf b}  \\
%{\bf g} \\
%\end{array} \right]
%$$
We need to use the following recurrence relations in an iterative process to produce a bidiagonal matrix
\begin{eqnarray}
\label{this}
AV_{k}&=&U_{k+1}B_{k} \\
\label{that}
A^{T}U_{k+1}&=&V_{k}B_{k}^{T}+\alpha_{k+1}{\bf v}_{k+1}{\bf e}_{k+1}^{T}
\end{eqnarray}
where $V_{k}$ and $U_{k}$ are matrices with orthonormal columns, and 
$$ B_{k} =\left[ \begin{array}{cccc}
\alpha_{1} & \, & \,  & \,  \\
\beta_{2} & \alpha_{2} & \, & \, \\
\, & \beta_{3} & \ddots & \, \\
\, & \, & \ddots & \alpha_{k} \\
\, & \, & \, & \beta_{k+1} \\
\end{array} \right].
$$
Also we have that 
\begin{eqnarray}
\nonumber
A^{T}AV_{k}&=&A^{T}U_{k+1}B_{k}=(V_{k}B_{k}^{T}+\alpha_{k+1}{\bf v}_{k+1}{\bf e}_{k+1}^{T})B_{k}\\
\label{2.9}
&=&V_{k}B_{k}^{T}B_{k}+\hat{\alpha}_{k}{\bf v}_{k+1}{\bf e}_{k+1}^{T}
\end{eqnarray}
and
$$
\hat{\alpha}_{k+1}=\alpha_{k+1}\beta_{k+1}.
$$

%$$V_{k}=[v_{1}\,v_{2}\, \ldots \v_{k}] \quad \quad $$U_{k}=[u_{1} \, u_{2} \, \ldots \,u_{k}]$$
Because $B_{k}$ is bidiagonal, it follows that $B_{k}^{T}B_{k}$ is symmetric and tridiagonal. It can be seen from (\ref{2.9}) that (\ref{this}) and (\ref{that})  implicitly apply Lanczos iteration to $A^{T}A$. 
Now this iterative process can be used to obtain the approximate solution to the forward and adjoint systems. We define the residuals at step $k$ as 
\begin{eqnarray}
\label{res1}
{\bf r}_{k}&=&{\bf b}-A{\bf x}_{k}\\
\label{res2}
{\bf s}_{k}&=&{\bf g}-A^{T}{\bf y}_{k}
\end{eqnarray}
 where
$${\bf x}_{k}={\bf x}_{0}+V_{k}{\bf z}_{k} \quad \quad {\bf y}_{k}={\bf y}_{0}+U_{k+1}{\bf w}_{k.}$$
The goal of the LSQR approach is to obtain an approximation that minimizes the norm of the residual. That is, the norm $\| {\bf r}_{k} \| =\| {\bf b}-A{\bf x}_{k} \|$ is minimized. When working with the forward and adjoint problems, this approach is limited due to the relationship between the starting vectors
$$
A^{T}{\bf u}_{1}=\alpha_{1}{\bf v}_{1} .
$$
The above relationship does not allow ${\bf v}_{1}$ to be chosen independently.  

\subsection{Generalized LSQR (GLSQR)}
The GSLQR method \cite{GSW,GLSQR} overcomes the disadvantages of the LSQR method by choosing starting vectors ${\bf u}_{1}=\frac{{\bf r}_{0}}{\|{\bf r}_{0}\|}$ and ${\bf v}_{1}=\frac{{\bf s}_{0}}{\|{\bf s}_{0}\|}$ independently where, for an initial guess of ${\bf x}_{0}$ and ${\bf y}_{0}$, 
${\bf r}_{0}={\bf b}-A{\bf x}_{0}$ and ${\bf s}_{0}={\bf g}-A^{T}{\bf y}_{0}.$
It is based on the factorizations
\begin{eqnarray}
AV_{k}&=&U_{k+1}T_{k+1,k}=U_{k}T_{k,k}+\beta_{k+1}{\bf u}_{k+1}{\bf e}_{k}^{T}\\
A^{T}U_{k}&=&V_{k+1}S_{k+1,k}=V_{k}S_{k,k}+\eta_{k+1}{\bf v}_{k+1}{\bf e}_{k}^{T}
\end{eqnarray}
From the above we get that
\begin{eqnarray}
\label{2.19}
\beta_{k+1}{\bf u}_{k+1}&=&A{\bf v}_{k}-\alpha_{k}{\bf u}_{k}-\gamma_{k-1}{\bf u}_{k-1}={\bf c}_{k}\\
\label{2.20}
\eta_{k+1}{\bf v}_{k+1}&=&A^{T}{\bf u}_{k}-\delta_{k}{\bf v}_{k}-\theta_{k-1}{\bf v}_{k-1}={\bf d}_{k},
\end{eqnarray}
where the recursion coefficients $\alpha_{k}$, $\gamma_{k}$, $\eta_{k}$, and $\theta_{k}$ are chosen to make $U_{k}$ and $V_{k}$ have orthonormal columns, which yields 
\begin{eqnarray}
\alpha_{k}&=&{\bf u}_{k}^{T}A{\bf v}_{k},\\
\gamma_{k}&=&{\bf u}_{k-1}^{T}A{\bf v}_{k+1}, \\
\delta_{k}&=&{\bf v}_{k}^{T}A^{T}{\bf u}_{k}, \\
\theta_{k}&=&{\bf v}^{T}A^{T}{\bf u}_{k+1}.
\end{eqnarray}
We can define ${\bf u}_{k+1}=\frac{{\bf c}_{k}}{\beta_{k}}$ and ${\bf v}_{k}=\frac{{\bf d}_{k}}{\eta_{k}}$, where $\beta_{k}=\|{\bf c}_{k}\|$, and $\eta_{k}=\|{\bf d}_{k}\|$.  
Now we have that 
$$ T_{k+1,k} =\left[ \begin{array}{cccc}
\alpha_{1} & \gamma_{1} & \,  & \,  \\
\beta_{2} & \alpha_{2} & \, & \, \\
\, & \ddots & \ddots & \gamma_{k-1} \\
\, & \, &  \beta_{k}& \alpha_{k} \\
\, & \, & \, & \beta_{k+1} \\
\end{array} \right]
\quad
 S_{k+1,k} =\left[ \begin{array}{cccc}
\delta_{1} & \theta_{1} & \,  & \,  \\
\eta_{2} & \delta_{2} & \ddots & \, \\
\, & \ddots & \ddots & \theta_{k-1} \\
\, & \, & \eta_{k} & \delta_{k} \\
\, & \, & \, & \eta_{k+1} \\
\end{array} \right].
$$
The residuals can be expressed as follows
\begin{equation}
\|{\bf r}_{k}\|=\|{\bf r}_{0}-U_{k+1}T_{k+1,k}{\bf x}_{k}\|=\| \|{\bf r}_{0}\|{\bf e}_{1}-T_{k+1,k}{\bf x}_{k}\|,
\end{equation}
and
\begin{equation}
\|{\bf s}_{k}\|=\|{\bf s}_{0}-V_{k}S_{k+1,K}^{T}{\bf y}_{k}-\alpha_{k+1}{\bf v}_{k+1}{\bf e}_{k+1}^{T}{\bf y}_{k}\|.
\end{equation}
The solutions ${\bf x}_{k}$ and ${\bf y}_{k}$ are 
\begin{eqnarray}
{\bf x}_{k}&=&{\bf x}_{0}+\|{\bf r}_0\|V_{k}T_{k,k}^{-1}{\bf e}_{1}\\
{\bf y}_{k}&=&{\bf y}_{0}+\|{\bf s}_{0}\|U_{k}S_{k,k}^{-1}{\bf e}_{1}.
\end{eqnarray}

\section{Iterative Methods for Nonsymmetric Saddle Point Matrices}
\label{sec4}
The matrix $M$, defined as follows
\begin{equation}
\label{M}
M\equiv 
\left[ \begin{array}{cc}
A^{T}WA & A^{T}  \\
-A & 0 \\
\end{array} \right], 
\end{equation} 
where $A \in \mathbb{R}^{n\times n}$ is invertible and $W$ is a symmetric positive definite matrix, is an example of a \textit{nonsymmetric saddle point matrix}. 
It can be shown that 
%if the matrix $W$ is symmetric positive definite, meaning that ${\bf y}^{T} W {\bf y} > 0$ for all ${\bf y} \neq$ {\bf 0}, then 
${\bf x}^{T} M {\bf x}\geq 0$ for all ${\bf x}\neq {\bf 0}$. To see this, we first let 
${\bf x}= \left[ \begin{array}{cc}
{\bf y}  \\
{\bf z} \\
\end{array} \right].$
Then 
${\bf x}^{T} M {\bf x}$ can be written as 
\begin{eqnarray*}
{\bf x}^{T}M{\bf x}& = &
\left[ \begin{array}{cc}
{\bf y}^{T} & {\bf z}^{T}  \\
\end{array} \right]
\left[ \begin{array}{cc}
A^{T}WA & A^{T}  \\
-A & 0 \\
\end{array} \right]
\left[ \begin{array}{cc}
{\bf y}  \\
{\bf z} \\
\end{array} \right] \\
& = & {\bf y}^{T}(A^{T}WA ) {\bf y}-{\bf z}^{T}A{\bf y}+{\bf y}^{T}A^{T}{\bf z} \\
& = & {\bf y}^{T}(A^{T}WA ) {\bf y}.
\end{eqnarray*}
%Carrying out this product, we get\,
%${\bf x}^{T}M{\bf x}=$. It is easy to see that this product (which results in a scalar value) reduces to \,${\bf y}^{T}(A^{T}WA ) {\bf y}$\, due to the cancellation of terms. 
%It is also worthwhile to note that, since $W$ is chosen to be symmetric positive definite, then %${\bf y}^{T}(A^{T}WA ) {\bf y}\, >0 $ when ${\bf y}$ is nonzero. \\

Now, if we let ${\bf r}=A{\bf y}$ for any nonzero vector ${\bf y}$, then, 
${\bf r}^{T}=(A{\bf y})^{T}={\bf y}^{T}A^{T}$. 
since $W$ is symmetric positive definite, we have that  ${\bf y}^{T}(A^{T}WA){\bf y}={\bf r}^{T}W{\bf r} >0$, since {\bf r} is nonzero due to $A$ being invertible. 
%Now that we have shown that ${\bf x}^{T}M{\bf x} \geq 0$, we need to show equality. To do this we note that while ${\bf x} \neq 0$, any component of {\bf x} could equal 0, so long as there is at least one non-zero component. 
On the other hand, if we assume ${\bf y} =0$, then ${\bf x}^{T}M{\bf x} ={\bf y}^{T}(A^{T}WA){\bf y}=0.$ That is, whether ${\bf y}$ is nonzero or not, ${\bf x}^{T}M{\bf x} ={\bf r}^{T}W{\bf r} \geq 0$.

\subsection{Ensuring a Real Positive Spectrum}

We want to choose $W$ so that the matrix $M$ has a real positive spectrum, 
%meaning it has all real, positive eigenvalues, 
so it is suitable for a conjugate gradient-like iteration \cite{LP}. To make this choice we need to first define 
$$
\mathcal{M}(\gamma) \equiv \mathcal{J}p(M) =\mathcal{J}(M-\gamma I)=
\left[ \begin{array}{cc}
A^{T}WA -\gamma I & A^{T} \\
A & \gamma I \\
\end{array} \right],
$$
where $p$ is a polynomial of degree one in the form $p(\zeta) = \zeta - \gamma$ for $ \gamma \in \mathbb{R}$ and 
$$
\mathcal{J} \equiv
\left[ \begin{array}{cc}
I & 0\\
0 & -I\\
\end{array} \right].
$$
The goal here is to determine if there exists a symmetric positive definite matrix $\mathcal{M}(\gamma)$ with respect to which $ M$ is symmetric, meaning that $M$ is $\mathcal{M}(\gamma)$-symmetric if $\mathcal{M}(\gamma)M=M^{T}\mathcal{M}(\gamma)=(\mathcal{M}(\gamma)M)^{T}$. 

Let us first define a generic nonsymmetric saddle point matrix 
$$\mathcal{A}=
\left[ \begin{array}{cc}
\hat{A} & \hat{B}^{T}\\
-\hat{B} & \hat{C} \\
\end{array} \right].$$ and then define $\mathcal{M}(\gamma)=\mathcal{J}p(\mathcal{A})$. 
We can use the following results from \cite{LP} to determine how to obtain a real positive spectrum:
\begin{lemma}
\label{lem2}
Let the matrix 
$$
\mathcal{J}\equiv 
\left[ \begin{array}{cc}
I & 0\\
0 & -I\\
\end{array} \right]
$$
be conformally partitioned with $\mathcal{A}$. Then\\
(1) $\mathcal{A}$ is $\mathcal{J}$-symmetric, i.e., $\mathcal{J}\mathcal{A}=\mathcal{A}^{T}\mathcal{J}=(\mathcal{J}\mathcal{A})^{T}$, and for any polynomial $p$,\\
(2) $p(\mathcal{A})$ is $\mathcal{J}$-symmetric, i.e., $\mathcal{J}p(\mathcal{A})=p(\mathcal{A}^{T})\mathcal{J}=(\mathcal{J}p(\mathcal{A}))^{T}$, and\\
(3) $\mathcal{A}$ is $\mathcal{J}p(\mathcal{A})$-symmetric, i.e., $(\mathcal{J}p(\mathcal{A}))\mathcal{A}=\mathcal{A}^{T}
(p(\mathcal{A})^{T})\mathcal{J}=(\mathcal{J}p(\mathcal{A})\mathcal{A})^{T}$.
\end{lemma}
\begin{theorem}
\label{thm1}
The symmetric matrix $\mathcal{M}(\gamma)$ is positive definite if and only if 
\begin{equation}
\label{4.1}
\lambda_{min}(\hat{A})>\gamma >\lambda_{max}(\hat{C})
\end{equation}
where $\lambda_{min}$ and $\lambda_{max}$ denote the smallest and largest eigenvalues, respectively, and 
\begin{equation}
\| (\gamma I-\hat{C})^{-1/2}\hat{B}(\hat{A}-\gamma I)^{-1/2} \|_{2} <1.
\end{equation}
\end{theorem}
A sufficient condition that makes $\mathcal{M}(\gamma)$ positive definite can be derived from the above theorem.
\begin{corollary}
\label{cor1}
The matrix $\mathcal{M}(\gamma)$ is symmetric positive definite when (\ref{4.1})
holds, and, in addition,
\begin{equation}
\label{4.4}
\|\hat{B}\|_{2} ^{2} < (\lambda_{min} (\hat{A})-\gamma )(\gamma-\lambda_{max}(\hat{C})).
\end{equation}
For $\gamma =\hat{\gamma}\equiv \frac{1}{2}(\lambda_{min}(\hat{A})+\lambda_{max}(\hat{C}))$, the right hand side of (\ref{4.4}) is maximal and (\ref{4.4}) reduces to 
\begin{equation}
\label{4.5}
2\|\hat{B}\|_{2} < (\lambda_{min}(\hat{A})-\lambda_{max}(\hat{C})).
\end{equation}
\end{corollary}
The preceding results lead to a simple approach to determining whether $\mathcal{A}$ is suitable for a conjugate gradient-like iteration \cite{LP}.
\begin{corollary}
\label{cor2}
If there exists a $\gamma \in \mathbb{R}$ so that $\mathcal{M}(\gamma)$ is positive definite, then $\mathcal{A}$ has a nonnegative real spectrum and a complete set of eigenvectors that are orthonormal with respect to the inner product defined by $\mathcal{M}(\gamma)$. In case $\hat{B}$ has full rank, the spectrum of $\mathcal{A}$ is real and positive. 
\end{corollary}
Using the previous results from \cite{LP}, we obtain a simple criterion for determining whether the matrix $M$ from (\ref{M}) can be constructed in such a way as to satisfy the criterion 
in Corollary \ref{cor2}.
\begin{theorem} \label{wthm} Let $A$ be an invertible $n\times n$ real matrix, %let the scalar $w$ satisfy
and let $W$ be a symmetric positive definite $n\times n$ matrix that satisfies
\begin{equation}
\label{Wcond}
\sigma_{\min}(W) > \frac{2\kappa_2(A)}{\sigma_{\min}(A)}.
\end{equation}
%and let $W=wI$.  
Then the matrix $M$ defined by
%If we assume $W$ to be symmetric positive definite a proper selection of $W$ can be made so that
$$M=
\left[ \begin{array}{cc}
A^{T}WA & A^{T} \\
-A & 0\\
\end{array} \right]
$$
has real positive eigenvalues and eigenvectors that are orthogonal with respect to the inner product defined by $\mathcal{M}(\gamma) = \mathcal{J}p(M)$. 
That is, the above selection of $W$ makes the matrix $M$ suitable for a conjugate gradient-like iteration.
\end{theorem}
Proof:
We need to satisfy (\ref{4.1}) with a proper selection of $\gamma$. Let 
$$\gamma=\frac{1}{2}(\lambda_{min}(A^{T}WA)),$$ based on Corollary \ref{cor1}. Because of how $\gamma$ is defined, $\gamma $ satisfies %the following
\begin{equation}
\lambda_{min}(A^{T}WA)>\gamma > 0,
\end{equation}
which means (\ref{4.1}) is also satisfied. 
Now we need to choose $W$ so that (\ref{4.5}) from Corollary \ref{cor1} holds. We require
\begin{equation}
2\|A^{T}\|_{2} < \lambda_{\min}(A^{T}WA),
\end{equation}
or
\begin{equation}
\label{4.8}
2\sigma_{\max}(A) <\lambda_{\min}(A^{T}WA)
\end{equation}
where $\|A^{T}\|_{2}=\|A\|_{2}$ is equal to the largest singular value of $A$, $\sigma_{\max}(A)$. 
From the fact that $A^T WA$ is symmetric positive definite, we obtain
\begin{eqnarray*}
\frac{1}{\lambda_{\min}(A^T WA)} & = & \rho( (A^T W A)^{-1} ) \\
& = & \|(A^T WA)^{-1}\|_2 \\
& \leq & \|A^{-1}\|_2^2 \|W^{-1}\|_2 \\
& \leq & \frac{1}{\sigma_{\min}(A)^2 \sigma_{\min}(W)}.
\end{eqnarray*}
Therefore, (\ref{4.8}) is satisfied if
\begin{equation}
2\sigma_{\max}(A) < \sigma_{\min}(A)^2 \sigma_{\min}(W),
\end{equation}
or, equivalently, if (\ref{Wcond}) is satisfied.
%If we let $W=wI$, then
%$$
%\lambda_{min}(A^{T}WA)=w\lambda_{\min}(A^{T}A).$$ 
%Rearranging (\ref{4.8}) gives the following criterion for $w$,
%\begin{equation}
%\label{w2}
%w> \frac{2\sigma_{\max}(A)}{\lambda_{\min}(A^{T}A)}.
%\end{equation}
$\Box$

It follows that the matrix $W$ satisfies the requirements to make $\mathcal{M}(\gamma)$ be symmetric positive definite and that $\mathcal{A}=M$ has a real, positive spectrum from Corollary \ref{cor2}. This result makes the matrix suitable for a conjugate gradient-like iteration, as will be described below.
% in the sense of Lemma \ref{lem1} below. 

\subsection{The Case $W=wI$}

Let $A=U\Sigma V^T$ be the SVD of $A$, where
$$U = \left[ \begin{array}{ccc} {\bf u}_1 & \cdots & {\bf u}_n \end{array} \right], \quad
V = \left[ \begin{array}{ccc} {\bf v}_1 & \cdots & {\bf v}_n \end{array} \right]$$
and $\Sigma = \mbox{diag}(\sigma_1,\ldots,\sigma_n)$.  In the case $W=wI$ for some scalar $w$, the condition from Theorem \ref{wthm} reduces to
\begin{equation}
\label{w}
w> \frac{2\kappa_2(A)}{\sigma_n}.
\end{equation}

We now study the eigensystem of $M$.  Let $M{\bf x}_j = \lambda_j {\bf x}_j$ for $j=1,2,\ldots,2n$, where
${\bf x}_j = \left[ \begin{array}{cc} {\bf y}_j^T & {\bf z}_j^T \end{array} \right]^T$.  The form of $M$ from (\ref{M}), with $W=wI$, yields
\begin{eqnarray}
wA^T A{\bf y}_j + A^T {\bf z}_j & = & \lambda_j {\bf y}_j, \label{eig1} \\
-A{\bf y}_j & = & \lambda_j {\bf z}_j \label{eig2}
\end{eqnarray}
for $j=1,2,\ldots,2n$.  Substituting (\ref{eig2}) into (\ref{eig1}) yields
\begin{equation}
(1-w\lambda_j) A^T {\bf z}_j = \lambda_j {\bf y}_j.
\end{equation}
Multiplying through by $A$ and applying (\ref{eig2}), we obtain
$$(w\lambda_j - 1) AA^T {\bf z}_j = \lambda_j^2 {\bf z}_j.$$
It follows that each ${\bf z}_j$ is a multiple of a left singular vector of $M$, and $\lambda_j^2 / (w\lambda_j - 1)$ is the square of the corresponding singular value.
Furthermore, from (\ref{eig2}), we find that ${\bf y}_j$ is a multiple of a right singular vector of $M$.

We conclude that the eigenvectors ${\bf x}_1,{\bf x}_2,\ldots, {\bf x}_{2n}$ of $M$ are given by
\begin{equation} \label{Mevcts}
{\bf x}_{2j-1} = \left[ \begin{array}{c} -\lambda_j^+ {\bf v}_j \\ \sigma_j {\bf u}_j \end{array} \right], \quad
{\bf x}_{2j} = \left[ \begin{array}{c} -\lambda_j^- {\bf v}_j \\ \sigma_j {\bf u}_j \end{array} \right], \quad j = 1, 2,\ldots, n,
\end{equation}
with corresponding eigenvalues $\lambda = \lambda_j^+, \lambda_j^-$ that satisfy the quadratic equation
\begin{equation} \label{Mevals}
\lambda^2 - \sigma_j^2 w \lambda + \sigma_j^2 = 0.
\end{equation}
It can be shown directly from (\ref{Mevcts}) and (\ref{Mevals}) that these eigenvalues are real and positive,
and the corresponding eigenvectors linearly independent, if and only if $w$ satisfies the weaker condition
\begin{equation}
w > \frac{2}{\sigma_n},
\end{equation}
which is consistent with the necessary and sufficient condition for $\mathcal{M}(\gamma)$ to be positive definite given in Theorem \ref{thm1}.
%Furthermore, from the eigenvalues 
%$$\lambda_j^+,\lambda_j^- = \frac{\sigma_j^2 w \pm \sigma_j\sqrt{\sigma_j^2 w^2 - 4}}{2},$$
%if we let $w = 2(1+\epsilon)/\sigma_n$, then we have
%$$\lambda_j^+,\lambda_j^- = \frac{1}{\sigma_n} \left[ \sigma_j^2 (1+\epsilon) \pm \sigma_j\sqrt{\sigma_j^2 (1+\epsilon)^2 - \sigma_n^2}\right],$$

\subsection{Nonsymmetric Saddle Point Conjugate Gradient Method}

Let $A\in\mathbb{R}^{n\times n}$ be nonsymmetric.  We will now introduce a Conjugate Gradient (CG) approach that solves the linear system $A{\bf x} = {\bf c}$
by solving an equivalent system of the form $M{\bf z} = {\bf b}$, where
\begin{equation}
M \equiv 
\left [\begin{array}{cc}
A^{T}W{A} &  A^{T}\\
-A & 0 \\
\end{array} \right].
\end{equation}
The matrix $M$ is also not symmetric; however, the spectrum is entirely contained in the right half of the complex plane, due to the fact that $ {\bf x}^{T}M{\bf x} \geq 0$ for all ${\bf x}$. 
In the preceding discussion, we established that if $W$ was chosen so as to satisfy the assumptions of Theorem \ref{wthm}, then
% We know that there exists a conjugate gradient-like method for solving systems with this matrix $M$ because 
$M$ is diagonalizable with real, positive eigenvalues.  Furthermore, the bilinear form $({\bf u},{\bf v})_{G}={\bf u}^{T}G{\bf v}$, where $G = \mathcal{M}(\gamma) = \mathcal{J}p(M)$, 
is a proper inner product, as $G$ is symmetric positive definite.  It follows that $M$ is $G$-symmetric and $G$-definite, meaning that $(M{\bf u},{\bf v})_{G}=({\bf u},M{\bf v}_{G})$
for all ${\bf u}, {\bf v}\in \mathbb{R}^{2n}$, and $({\bf u},M{\bf u})_G > 0$ for all ${\bf u}\neq {\bf 0}$. 

Let the vectors ${\bf p}$ and ${\bf b}$ be defined by 
\begin{equation}
 {\bf b}=
\left[\begin{array}{c}
A^{T}W{\bf c}+{\bf d}\\
-{\bf c}\\
\end{array} \right],
\quad 
{\bf p}=
\left[ \begin{array}{c}
{\bf d} \\
{\bf 0} \\
\end{array} \right],
\end{equation}
where $A{\bf x}={\bf c}$, $M{\bf z} = {\bf b}$, and ${\bf p}^{T}{\bf z} = {\bf d}^T {\bf x}$ 
is the scattering amplitude for given vectors ${\bf c}$ and ${\bf d}$ that represent the field and antenna, respectively. The following conjugate gradient method is based on a given inner product 
$(\textbf{u},\textbf{v})_{G}=\textbf{v}^{T}G\textbf{u}$ for solving the linear system $M{\bf x}={\bf b}$. 
\begin{algorithmic}
\State {${\bf Algorithm \,\,3.1}$}
\State {{\bf Input:} System matrix $M$, right hand side vector ${\bf b}$, inner product matrix $W$, initial guess ${\bf x}_{0}$}
\vspace{.1cm}
\State {\bf Require:} {${\bf r}_{0}={\bf b}-M{\bf x}_{0}$}
\For {$i=0,1,\ldots$ until convergence}
\State {$\alpha_{i}=\frac{({\bf x}-{\bf x}_{i},{\bf p}_{i})_{G}}{({\bf p}_{i},{\bf p}_{i})_{G}}$}
\State {${\bf x}_{i+1}={\bf x}_{i}+\alpha_{i}{\bf p}_{i}$}
\State {${\bf r}_{i+1}={\bf r}_{i}-\alpha_{i}M{\bf p}_{i}$}
\State {$\beta_{i+1}=-\frac{({\bf r}_{i+1},{\bf p}_{i})_{G}}{({\bf p}_{i},{\bf p}_{i})_{G}}$}
\State {${\bf p}_{i+1}={\bf r}_{i+1}+\beta_{i+1}{\bf p}_{i}$}
\EndFor
\end{algorithmic}
We have the inner product matrix $G=\mathcal{M}(\gamma)M$ suggested by \cite{LP}. From \cite{LP}, we see that this choice of $G$ gives a working CG from the following lemma.
\begin{lemma}
\label{lem1}
Suppose that the symmetric matrix $\mathcal{M}(\gamma)$ is positive definite. Then Algorithm 3.1 is well defined for $M$ and $G=\mathcal{M}(\gamma)M$, and (until convergence) the scalars $\alpha_{i}$ and $\beta_{i+1}$ can be computed as 
\begin{equation}
\label{alpha}
\alpha_{i}=\frac{({\bf r}_{i}, {\bf r}_{i})_{\mathcal{M}(\gamma)}}{(M{\bf p}_{i},{\bf p}_{i})_{\mathcal{M}(\gamma)}}
\end{equation}
\begin{equation}
\beta_{i+1}=\frac{({\bf r}_{i+1}, {\bf r}_{i+1})_{\mathcal{M}(\gamma)}}{(M{\bf r}_{i},{\bf r}_{i})_{\mathcal{M}(\gamma)}}.
\end{equation}
\end{lemma}
%Since we are using nonsymmetric saddle point conjugate gradient method, we need each residual to be orthogonal to each previous residual.
With this choice of inner product matrix, it can be shown that the residuals computed using the preceding algorithm are, in some sense, orthogonal.
\begin{theorem}
Each residual ${\bf r}_{k}$ as defined in Algorithm 3.1 is orthogonal to all previous residuals with respect to $\mathcal{M}(\gamma)$ ,i.e. $({\bf r}_{i}^{T},{\bf r}_{j})_{\mathcal{M}(\gamma)}=0$, where $i \neq j$.
\end{theorem}
Proof: We know that $ {\bf r}_{i+1}={\bf r}_{i}-\alpha_{i}M{\bf p}_{i}$. Let $\alpha_{i}$ be defined as in (\ref{alpha}). Also, we know that all of the search directions are orthogonal, i.e. ${\bf p}_{i}^{T}\mathcal{M}(\gamma)M{\bf p}_{j}=0$ for $i \neq j$. We want to show that ${\bf r}_{i}\mathcal{M}(\gamma){\bf r}_j =0$. This will be shown by induction, where the base case that we need to establish is 
\begin{equation}
\label{base}
{\bf r}_{i+1}^{T}\mathcal{M}(\gamma){\bf r}_{i}=0, \quad i=0,1, \ldots.
\end{equation}
To show this we use the definition of $\alpha_{i}$ and the expression for the search directions in the above algorithm, $ {\bf r}_{i+1}={\bf r}_{i}-\alpha_{i}M{\bf p}_{i}$. 
Now we have that 
\begin{equation}
\label{step1}
{\bf r}_{i+1}^{T}\mathcal{M}(\gamma){\bf r}_{i}={\bf r}_{i}^{T}\mathcal{M}(\gamma){\bf r}_{i}-\frac{{\bf r}_{i}^{T}\mathcal{M}(\gamma){\bf r}_{i}}{{\bf p}_{i}^{T}M^{T}\mathcal{M}(\gamma){\bf p}_{i}}{\bf p}_{i}^{T}M^{T}\mathcal{M}(\gamma){\bf r}_{i}.
\end{equation}
Reindexing the definition of the residual from the algorithm yields the following expression for ${\bf r}_{i}$
$$
{\bf r}_{i}={\bf p}_{i}-\beta_{i}{\bf p}_{i-1}.
$$
Substituting this into (\ref{step1}) gives
\begin{equation}
\label{step2}
{\bf r}_{i+1}^{T}\mathcal{M}(\gamma){\bf r}_{i}={\bf r}_{i}^{T}\mathcal{M}(\gamma){\bf r}_{i}-\frac{{\bf r}_{i}^{T}\mathcal{M}(\gamma){\bf r}_{i}}{{\bf p}_{i}^{T}M^{T}\mathcal{M}(\gamma){\bf p}_{i}}({\bf p}_{i}^{T}M^{T}\mathcal{M}(\gamma){\bf p}_{i}-\beta_{i}{\bf p}_{i-1}^{T}M^{T}\mathcal{M}(\gamma){\bf p}_{i})
\end{equation}
Rearranging the last term in (\ref{step2}) yields
$$
{\bf p}_{i-1}^{T}M^{T}\mathcal{M}(\gamma)\beta_{i}{\bf p}_{i}=\beta_{i}{\bf p}_{i}^{T}\mathcal{M}(\gamma)M{\bf p}_{i-1}=0
$$
because $\mathcal{M}(\gamma)$ is symmetric, and we already know that the search directions ${\bf p}_{i}$ are orthogonal with respect to $\mathcal{M}(\gamma)$. Now it is easy to see that the denominator in (\ref{step2}) and the last factor in the numerator cancel leaving
$$
{\bf r}_{i+1}^{T}\mathcal{M}(\gamma){\bf r}_{i}={\bf r}_{i}^{T}\mathcal{M}(\gamma){\bf r}_{i}-{\bf r}_{i}^{T}\mathcal{M}(\gamma){\bf r}_{i}=0.
$$
Now we need to show that each residual is orthogonal to all previous residuals. We will do this by showing ${\bf r}_{i}^{T}\mathcal{M}(\gamma){\bf r}_{i-d}=0$, where $d >1$. Our induction hypothesis is ${\bf r}_{i-1}^{T}\mathcal{M}(\gamma){\bf r}_{i-d}=0.$ To show this, first shift the indices to get the expression $${\bf r}_{i}={\bf r}_{i-1}-\alpha_{i-1}M{\bf p}_{i-1}.$$ Rearranging the recurrence relation for the search directions yields $${\bf r}_{i-d}={\bf p}_{i-d}-{\bf p}_{i-1-d}\beta_{i-d}.$$
Using this expression for ${\bf r}_{i}$ and ${\bf r}_{i-d}$ we obtain
\begin{eqnarray}
\nonumber
{\bf r}_{i}^{T}\mathcal{M}(\gamma){\bf r}_{i-d}&=&{\bf r}_{i-1}^{T}\mathcal{M}(\gamma){\bf r}_{i-d}-\alpha_{i-1}{\bf p}_{i-1}^{T}M\mathcal{M}(\gamma)({\bf p}_{i-d}-{\bf p}_{i-1-d}\beta_{i-d})\\
\nonumber
&=&{\bf r}_{i-1}^{T}\mathcal{M}(\gamma){\bf r}_{i-d}-\alpha_{i-1}{\bf p}_{i-1}^{T}M^{T}\mathcal{M}(\gamma){\bf p}_{i-d}+\\
& & \alpha_{i-1}{\bf p}_{i-1}^{T}M^{T}\mathcal{M}(\gamma){\bf p}_{i-1-d}\beta_{i-d},
\end{eqnarray}
where $${\bf r}_{i-1}^{T}\mathcal{M}(\gamma){\bf r}_{i-d}=0$$ by the induction hypothesis. 
Now we are left with 
$$
{\bf r}_{i}^{T}\mathcal{M}(\gamma){\bf r}_{i-d}=-\alpha_{i-1}{\bf p}_{i-1}^{T}M^{T}\mathcal{M}(\gamma){\bf p}_{i-d}+\alpha_{i-1}{\bf p}_{i-1}^{T}M^{T}\mathcal{M}(\gamma){\bf p}_{i-1-d}\beta_{i-d}=0,
$$
where both terms are 0 due to the orthogonality of the search directions.
$\Box$

\section{Numerical Results}

%In example \ref{ex1ilu} NspCG converges in about 70 iterations. GLSQR still fails to converge as before, and QMR takes about 200 iterations to reach an acceptable level of accuracy. 

\label{secna}
In this section, we will analyze the results from the methods described in this paper. These methods include QMR from Section 2.2, GLSQR from Section 2.3, and NspCG from Section 3.1. We have duplicated the results from \cite{GSW} for GLSQR and QMR and will compare them against the results for our NspCG method.
%We then show results for symmetrized initial vectors with $C$, unsymmetric Lanczos with $C$ and $M$ from Section 4.2 (bilinear form), Block GLSQR with $C$ from Section 3.4, perturbation with $M$ from Section 4.2, and unsymmetric block Lanczos with $M$ from Section 4.3. This chapter also includes an analysis of the eigenvalues of $M$ and $C$. 
%We can get a better idea of why particular methods perform in certain ways by looking at the eigenvalues of those matrices.

We need to first define the following matrix, 
%$$ 
%C=\left[ \begin{array}{cc}
%0 & A  \\
%A^{T}& 0\\
%\end{array} \right]
%$$
%where $C$ is a symmetric matrix, and 
$M$ is our nonsymmetric saddle point matrix
$$
M=\left[\begin{array}{cc}
A^{T}WA & -A \\
A^{T} & 0 \\
\end{array} \right]
$$
where $W=wI$ is defined from (\ref{w}). These examples are from \cite{GSW}. 
\subsection{Example 1}
This example uses the matrix created by
{\tt A=sprand(n,n,0.2)+speye(n)}
in {\sc Matlab} where \texttt{n=100}.
%, and the maximum number of iterations is 200. 
This creates a random sparse $n\times n$ matrix, where 0.2 is the density of uniformly distributed nonzero entries, and adds this to the identity. 
%\newpage
\begin{figure}[ht]
\centering
\includegraphics[width=10cm]{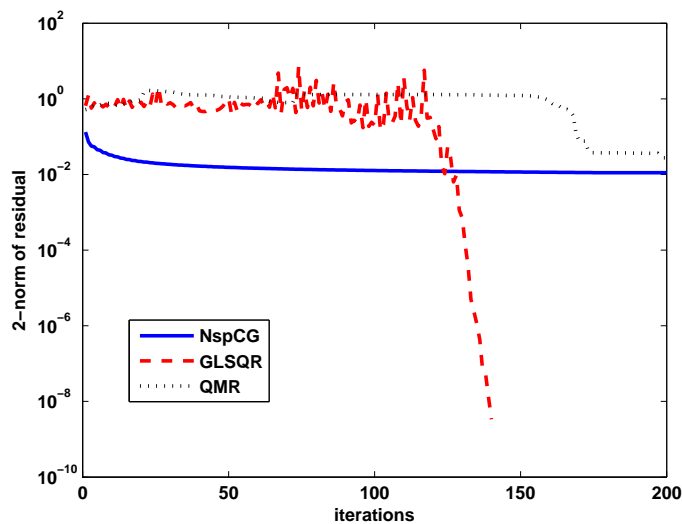}
\caption{Example 1 with the matrix $A$}
\label{ex1}
\end{figure}

In  Figure \ref{ex1} we see that at the beginning of the iteration NspCG reaches a better approximation in fewer iterations than either QMR or GLSQR. Although GLSQR eventually outperforms NspCG, it takes about 120 iterations before it shows any sign of convergence at all. Then it converges rapidly. 
%
%\begin{figure}[ht]
%\centering
%\includegraphics[width=12cm]{example1wC.eps}
%\caption{Example 1 with the matrix $C$}
%\label{ex1C}
%\end{figure}
%
%%\newpage
%We can see that Block GLSQR outperforms both the symmetrized method with the symmetric matrix $C$, and bilinear form with $C$. However, bilinear form $C$ does much better than symmetrized with $C$ reaching the maximum number of iterations without showing any sign of convergence. 

%\begin{figure}[ht]
%\centering
%\includegraphics[width=10cm]{Example1wM-eps-converted-to.pdf}
%\caption{Example 1 with the matrix $M$}
%\label{ex1M}
%\end{figure}

%% \newpage
%Figure \ref{ex1M} shows that bilinear form with $M$ oscillates a little bit but starts to converge toward 200 iterations. Perturpation with $M$ appears to do something similar, but it has several more sudden spikes throughout the iterations. Unsymmetric block Lanczos with $M$ appears to fail for this example. 
\subsection{Example 2}

Example 2 uses the {\sc ORSIRR\_1 } matrix from the Matrix Market collection, which represents a linear system used in oil reservoir modeling. This matrix can be obtained from {\tt http://math.nist.gov/MatrixMarket/}.
\begin{figure}[ht]
\centering
\includegraphics[width=10cm]{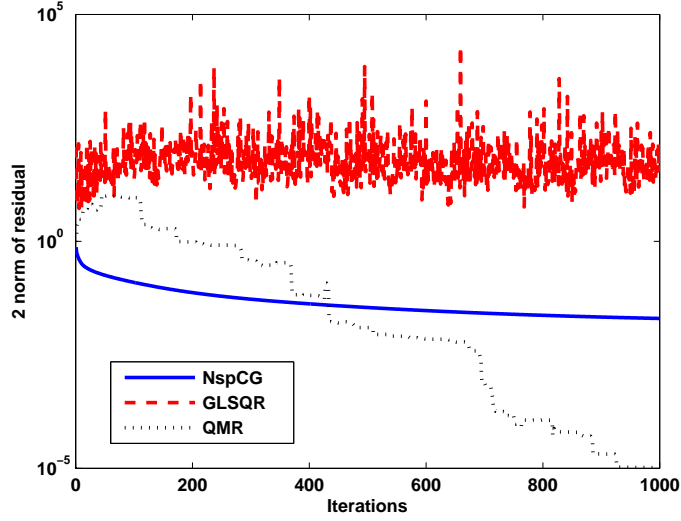}
\caption{Example 2}
\label{ex2}
\end{figure}

%\newpage

We see that NspCG starts out with the lowest error in the 2-norm of the residual. Also we see that in both Figure \ref{ex2} and Figure \ref{ex1} that NspCG is more consistent than either GLSQR or QMR. Although QMR actually outperforms GLSQR and NspCG, it takes about 400 iterations to do so.

%\begin{figure}[hb]
%\centering
%\includegraphics[width=12cm]{Example2wC-eps-converted-to.pdf}
%\caption{Example 2 with the matrix $C$}
%\label{ex2wC}
%\end{figure}
%
%From Figure \ref{ex2wC} it can be seen that bilinear form with $C$ is the only method showing any sign of convergence. Both block GLSQR with and symmetrized with $C$ fail in this case.
%
%\begin{figure}[ht]
%\centering
%\includegraphics[width=10cm]{Example2wM-eps-converted-to.pdf}
%\caption{Example 2 with the matrix $M$}
%\label{ex2wM}
%\end{figure}
%
%We see in Figure \ref{ex2wM} that perturbation with $M$ is the only method that converges, although it takes about 200 iterations to do so. Both bilinear form with $M$ and unsymmetric Lanczos with $M$ fail. 

\subsection{Example 3}

First define the circulant matrix
$$
J=
\left[ \begin{array}{cccc}
0 & 1 & \, & \, \\
\, & 0 & \ddots & \, \\
\, & \, & \ddots & 1 \\
1 & \, & \, & 0 \\
\end{array} \right].
$$
Now the matrix used in this example {\tt A=1e-3*sprand(n,n,0.2)+J}, where {\tt n=100}, can be constructed in {\sc Matlab}. 
%\newpage
\begin{figure}[ht]
\centering
\includegraphics[width=10cm]{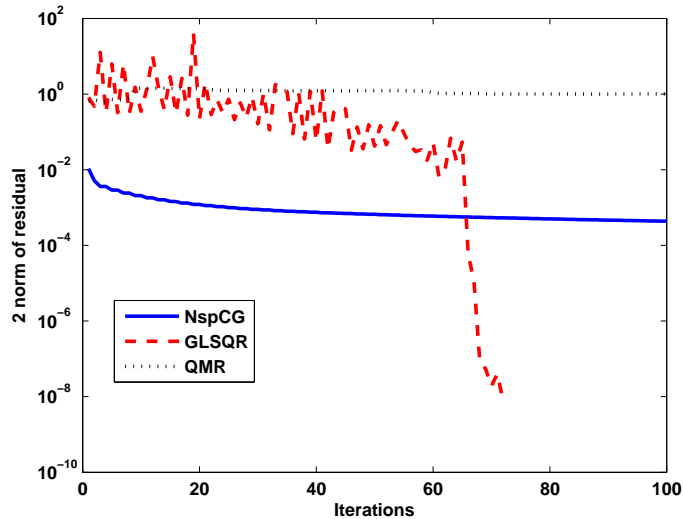}
\caption{Example 3 with the matrix $A$}
\label{ex3}
\end{figure}

NspCG starts out steady and consistent again in this Figure \ref{ex3} as we see in Figure \ref{ex2} and Figure \ref{ex1}. Eventually, GLSQR converges, taking about 70 iterations to do so, while QMR fails to show any sign of convergence.

%\begin{figure}[ht]
%\centering
%\includegraphics[width=10cm]{Example3wC-eps-converted-to.pdf}
%\caption{Example 3 with matrix $C$}
%\label{ex3C}
%\end{figure}
%
%In this Figure \ref{ex3C}, block GLSQR outperforms symmetrized with C and bilinear form with C by converging at about 140 iteration. Bilinear form with $C$ also shows some signs of convergence, but will take more than 200 iterations to converge, and the error in symmetrized with $C$ seems to be slowly descreasing. 
%
%\begin{figure}[ht]
%\centering
%\includegraphics[width=10cm]{Example3wM-eps-converted-to.pdf}
%\caption{Example 3 with matrix $M$}
%\label{ex3M}
%\end{figure}
%
%Figure \ref{ex3M} shows that unsymmetric block Lanczos with $M$, perturbation with $M$, and bilinear form with $M$ converge rapidly.  

\subsection{Example 4}

We need to first define 
$$
D_{1}=
\left[ \begin{array}{cccc}
1000 & \, & \, \\
\, & \ddots & \, \\
\, & \, & 1000 \\
\end{array} \right] \in \mathbb{R}^{p,p}
\quad \quad
D_{2}=
\left[ \begin{array}{cccc}
1 & \, & \,& \, \\
\, & 2 & \, & \, \\
\, & \, & \ddots & \, \\
\, & \, & \, & q \\
\end{array} \right] \in \mathbb{R}^{q,q}
$$
where $n=p+q$ and $\Sigma =$ diag$(D_{1}, D_{2})$. Now we can define $A=U\Sigma V^{T}$, where $U$ and $V$ are orthogonal matrices. For this example we use $n=100$ and $D_{1} \in \mathbb{R}^{90,90}$.
%\newpage
\begin{figure}[ht]
\centering
\includegraphics[width=10cm]{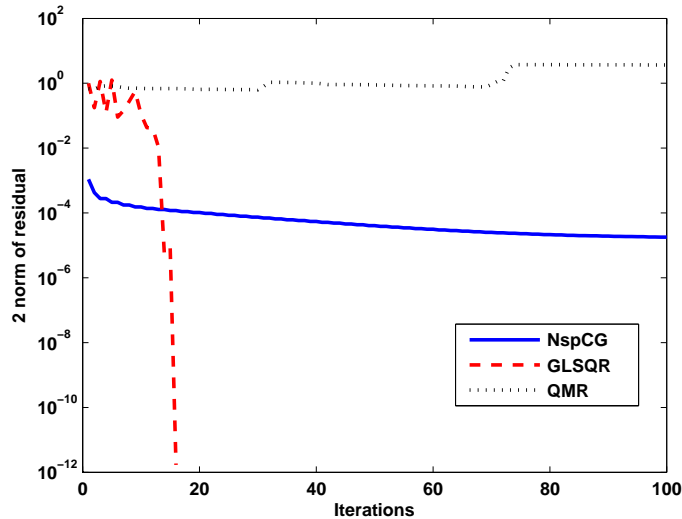}
\caption{Example 4 with the matrix $A$}
\label{ex4}
\end{figure}
From \ref{ex4} we see that NspCG starts off with the best approximation, but only for about 15 iterations. Then it is overtaken by GLSQR. Also, we can see that QMR fails to converge at all.

%\begin{figure}[ht]
%\centering
%\includegraphics[width=10cm]{Example4wC-eps-converted-to.pdf}
%\caption{Example 4 with the matrix $C$}
%\label{ex4C}
%\end{figure}
%
%In Figure \ref{ex4C}
%we see Block GLSQR and symmetrized with C both converged rapidly while it takes symmetrized with $C$ a few more iterations to do so. This is unlike the results from previous examples where we saw that bilinear form with $C$ outperformed symmetrized with $C$. This is likely due to the spreading out of eigenvalues of the matrix $C$. In this example, symmetrized with $C$ works better because most of the eigenvalues are clustered around 1000, whereas in other examples the eigenvalues are more spread out. 

%\begin{figure}[ht]
%\centering
%\includegraphics[width=10cm]{Example4wM-eps-converted-to.pdf}
%\caption{Example 4 with the matrix $M$}
%\label{ex4M}
%\end{figure}
%Figure \ref{ex4M} shows that perturbation with $M$ and bilinear form with $M$ converge rapidly, while unsymmetric Lanczos with $M$ fails to converge.
%%\newpage
\subsection{Example 5}
This example uses the same definition of $D_{1}$, $D_{2}$, and $A$ from Example 4. In this example we will let $n=100$ again, and $D_{1} \in \mathbb{R}^{50,50}$.
\begin{figure}[ht]
\centering
\includegraphics[width=10cm]{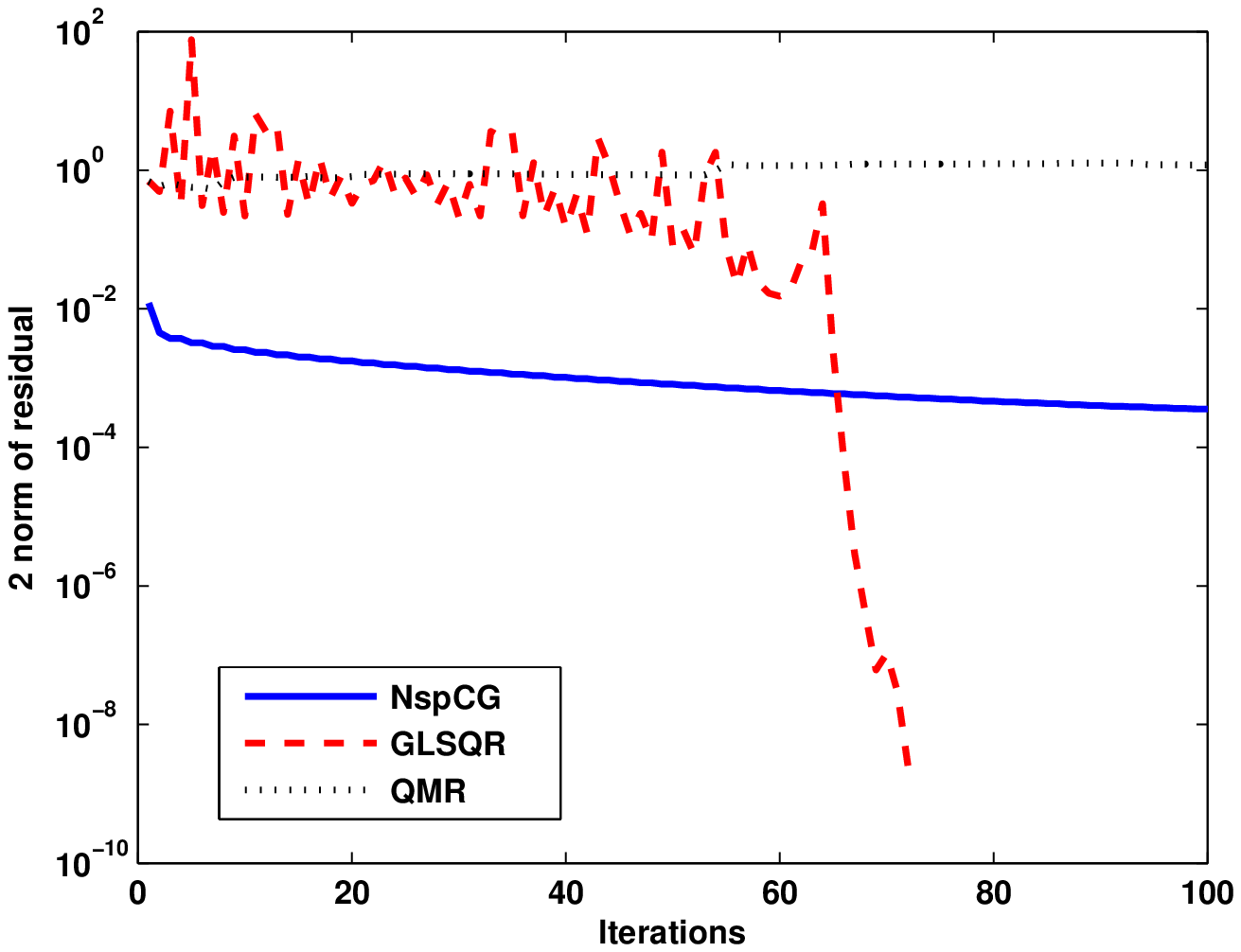}
\caption{Example 5 with the matrix $A$}
\label{ex5}
\end{figure}
Figure \ref{ex5} shows the same trend we have been seeing, that NspCG is more consistent at the beginning than any other method. At about 65 iterations GLSQR outperforms NspCG, and QMR fails to converge again.
%%\newpage
%\begin{figure}[ht]
%\centering
%\includegraphics[width=12cm]{Example5withC-eps-converted-to.pdf}
%\caption{Example 5 with matrix C}
%\label{ex5C}
%\end{figure}
%
%It can be seen from Figure \ref{ex5C} that block GLSQR with $C$ outperforms both methods. Symmetrized with $C$ does eventually converge at about 150 iterations, while bilinear form with $C$ fails to converge at all. 

%\begin{figure}[ht]
%\centering
%\includegraphics[width=10cm]{Example5withM-eps-converted-to.pdf}
%\caption{Example 5 with the matrix $M$}
%\label{ex5M}
%\end{figure}
%From Figure \ref{ex5M} we see that it takes more iterations in this example, but bilinear form with $M$, and perturbation with $M$ eventually start to converge around 100 iteration. In this example unsymmetric Lanczos with $M$ fails to converge. 

\subsection{Example 6}
This example uses the same definition of $D_{1}$, $D_{2}$, and $A$ from Example 4. In this example we will let $n=1000$ again, and $D_{1} \in \mathbb{R}^{600,600}$.
\begin{figure}[ht]
\centering
\includegraphics[width=10cm]{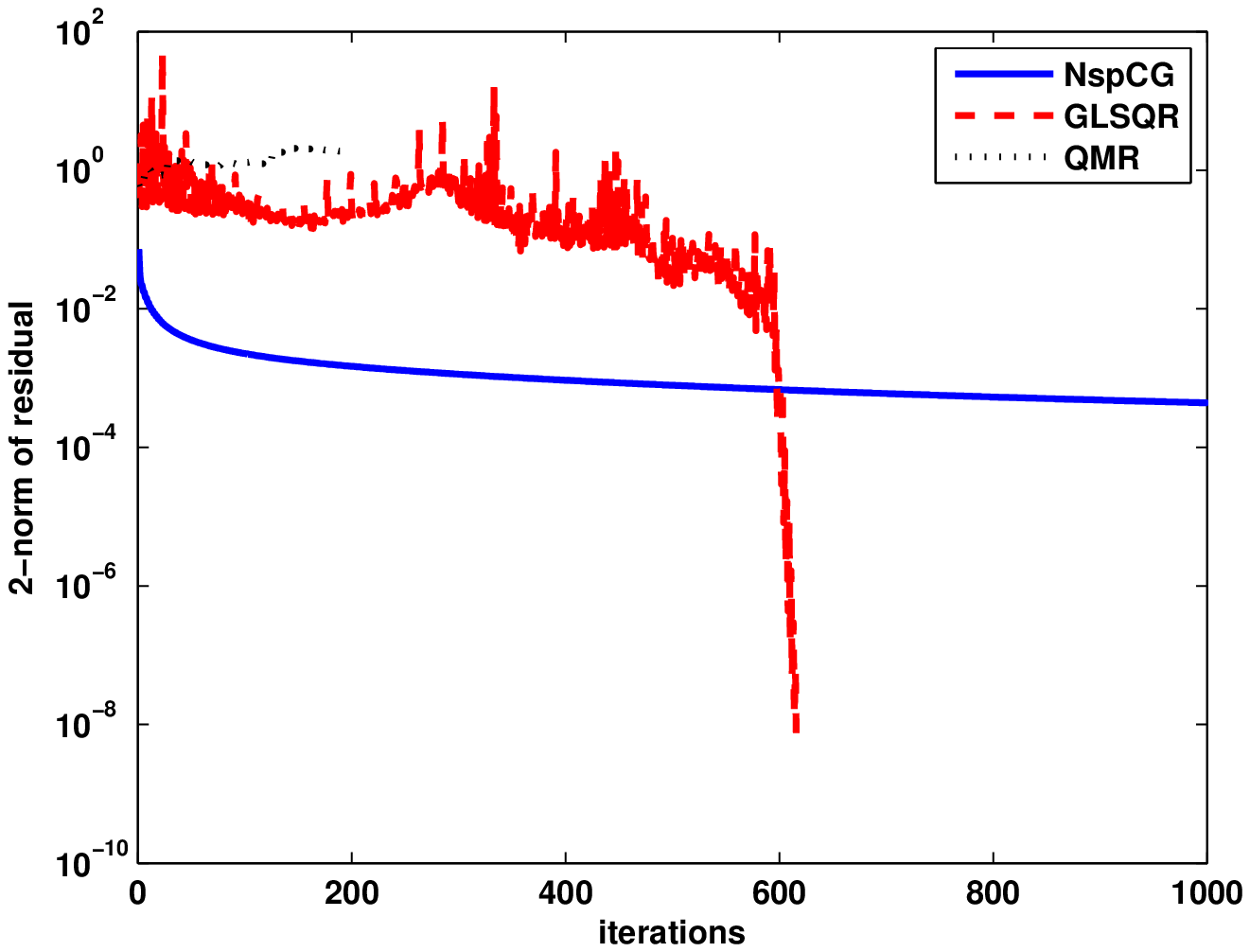}
\caption{Example 6 with the matrix $A$}
\label{ex6}
\end{figure}
From Figure \ref{ex6} we see that NspCG shows the best results for the first 600 iterations. GLSQR takes many iterations to converge in this case, and QMR does not converge at all.

%\subsection{Analysis of Eigenvalues of the Matrices $M$ and $C$}
%\begin{figure}[ht]
%\centering
%\includegraphics[width=10cm]{eigs1-eps-converted-to.pdf}
%\caption{Eigenvalues of $M$ and $C$ from Example 1}
%\label{ex1eigs}
%\end{figure}
%
%In Figure \ref{ex1eigs} the eigenvalues of $M$ from Example 1 are much more spread out than the eigenvalues of $C$ from Example 1. 
%
%\begin{figure}[ht]
%\centering
%\includegraphics[width=10cm]{eigs1pre-eps-converted-to.pdf}
%\caption{Eigenvalues of $M$ and $C$ from Example 1 with sparse matrix multiplied by $0.1$.}
%\label{ex1eigspre}
%\end{figure}
%We simulate the effect of preconditioning, as done in Figure \ref{ex1eigspre}, by transforming $A$ so that the eigenvalues of $M$ and $C$ are compressed. Now we can look at Example 1 with the random sparse matrix multiplied by $0.1$, i.e.,
%the matrix created by
%{\tt A=0.1*sprand(n,n,0.2)+speye(n)} in {\sc Matlab}, we get Figure \ref{Example1pre}.
%\begin{figure}[ht]
%\centering
%\includegraphics[width=12cm]{Simulatingpre-eps-converted-to.pdf}
%\caption{Example 1 with $M$ (eigenvalues compressed)}
%\label{Example1pre}
%\end{figure}
%From Figure \ref{Example1pre} we see that when the eigenvalues are compressed, modified CG converges rapidly. 

\section{Preconditioning}
\label{sec5b}
According to \cite{GVL} conjugate gradient has very rapid convergence for a symmetric positive definite matrix $A$ that is nearly identity. We need to apply {\em preconditioning} techniques to make our matrix $M$ satisfy this criterion. The result will be that the original system is transformed into an equivalent system where the coefficient matrix is near identity. As we have seen previously with conjugate gradient, preconditioning techniques can be generalized to the nonsymmetric case. The goal is to apply $ ILU$ preconditioning \cite{Saad}, while taking into account the structure of the nonsymmetric saddle point matrix $M defined in (\ref{M})$. The matrix $W$ in the (1,1) block is assumed to be a symmetric positive definite matrix; therefore it has a Cholesky factorization $W=GG^{T}$. We can use the $QR$ factorization $$G^{T}A=	QR$$ to obtain the factorization $M=LU$, where
\begin{equation}
L=
\left[ \begin{array}{cc}
R^{T} & 0 \\
-G^{-T}Q & G^{-T}Q 
\end{array} \right], \quad
U=
\left[ \begin{array}{cc}
R & Q^{T}G^{-1} \\
0 & Q^{T}G^{-1} 	
\end{array} \right].
\end{equation}
Let us define 
$$C= G^{T}A\widetilde{R}^{-1} \approx \widetilde{Q},$$ where an incomplete $QR$ factorization \cite{ILU} is computed from the sparse matrix $G^{T}A$ which gives $G^{T}A \approx \widetilde{Q}\widetilde{R}$ . By finding 
\begin{equation}
\widetilde{L}^{-1}=
\left[ \begin{array}{cc} 
\widetilde{R}^{-T} & 0 \\
\widetilde{R}^{-T} & \widetilde{Q}G^{T}
\end{array} \right], \quad
\widetilde{U}^{-1}=
\left[ \begin{array}{cc}
\widetilde{R}^{-1} & -\widetilde{R}^{-1} \\
0 & G\widetilde{Q}^{-T} 
\end{array} \right]
\end{equation}
it can be seen that the resulting preconditioned system matrix is given by
\begin{equation}
\widetilde{L}^{-1}M\widetilde{U}^{-1}=
\left[ \begin{array}{cc}
C^{T}C & -C^{T}C+C^{T}\widetilde{Q}\\
C^{T}C-\widetilde{Q}^{T}C & -C^{T}C+C^{T}\widetilde{Q}+\widetilde{Q}^{T}C
\end{array} \right].
\end{equation} 
The above matrix has the structure similar to that of $M$ from (\ref{M}), therefore it is a nonsymmetric saddle point matrix that is near $I$. 
\subsection{Example 1}
The following is Example 1 from the previous section with preconditioning.
\begin{figure}[ht]
\centering
\includegraphics[width=10cm]{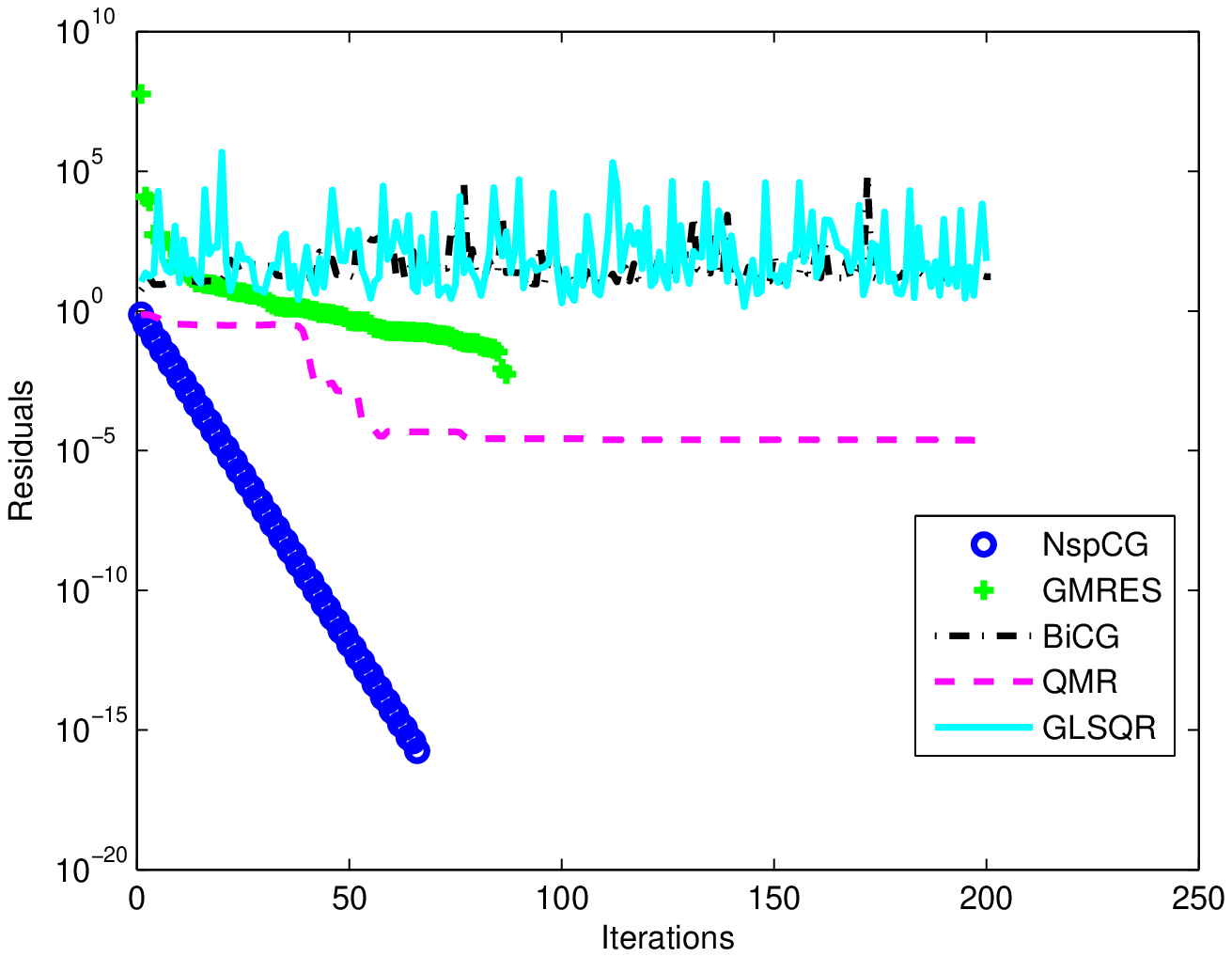}
\caption{Example 1 with preconditioning}
\label{ex1preconditioning}
\end{figure}
\subsection{Example 2}
The following is Example 2 from the previous section with preconditioning.
\begin{figure}[ht]
\centering
\includegraphics[width=10cm]{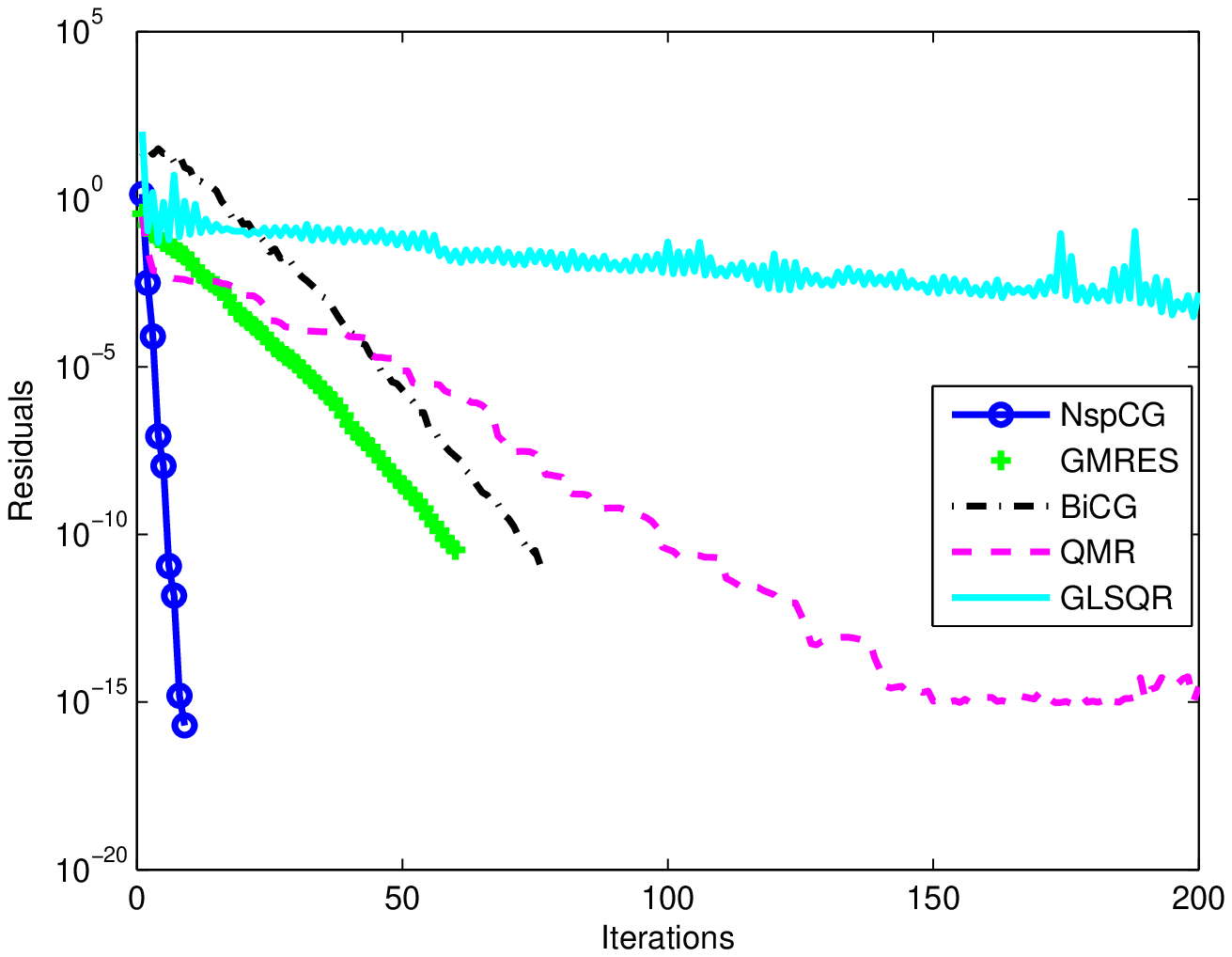}
\caption{Example 2 with preconditioning}
\label{ex2preconditioning}
\end{figure}
In Figure \ref{ex2preconditioning} NspCG converges very rapidly in only 10 iterations. QMR takes over 200 iterations, but still doesn't reach the level of accuracy that NspCG achieves. GLSQR does not converge at all. 

\section{Conclusions and Future Work}
\label{sec6}

The results from this paper show that the NspCG method is much more consistent and reliable than GLSQR or QMR. NspCG only takes a few iterations to make fairly significant progress while GLSQR takes many iterations in most cases, and QMR rarely makes any progress. 
If preconditioning is used with NspCG, as is usually done with a conjugate gradient method, we have provided evidence that it will dramatically accelerate convergence, compared to state-of-the-art
iterative methods such as GMRES or BiCG that are typically used to solve such systems.  These results support our hypothesis that more rapid convergence can be achieved by solving a system that,
while still nonsymmetric, shares essential properties with symmetric positive definite matrices and therefore is more suitable for conjugate gradient-like iteration.
%This is done by making the matrix closer to the identity. In Example 1 this is done by multiplying the random sparse matrix that is added to $I$ by a small factor, like 0.1, thus simulating the effect of preconditioning. A future goal would be to combine modified CG with preconditioning.

Future work will include relating the NspCG method to a quadrature rule, as in \cite{GM,MMQ4}, that can be used to compute the scattering amplitude without explicitly solving the forward or adjoint problem.  This has been done in \cite{GSW} with the symmetric matrix
$$C = \left[ \begin{array}{cc} 0 & A^T \\ A & 0 \end{array} \right]$$ in conjunction with block Lanczos iteration \cite{BLOCK}, but our goal is to achieve more rapid convergence.
Furthermore, because the forward system $A{\bf x} = {\bf b}$ is replaced with a system with twice as many unknowns and equations, it is essential to implement the iteration carefully
so that the gain in convergence speed is not offset by the additional expense of each iteration.  To that end, it is worthwhile to consider other choices for the matrix $W$ instead of just a multiple of identity.
%There are some difficulties in trying to do this because we have orthogonality with respect to the inner product matrix $W$ in this case. However we wish to try relate NspCG to a quadrature rule that would compute the scattering amplitude directly. The scattering amplitude is defined in terms of the standard inner product; it is this incompatibility of inner products that makes relation of NspCG to a quadrature rule difficult. 


%% \subsection{}

%% References


\bibliographystyle{elsarticle-num}



%\bibliographystyle{plain}

\end{document}